\documentclass[reqno, a4paper, 12pt]{amsart}
\usepackage{amsmath, amssymb, eucal, amscd, amstext, enumerate}
\usepackage{amsfonts,amsthm}
\usepackage{mathrsfs}

\usepackage[numbers,sort&compress]{natbib}
\usepackage{color,hyperref}

\theoremstyle{plain}
\newtheorem{thm}{Theorem}[section]
\newtheorem{lem}[thm]{Lemma}

\newtheorem{defn}[thm]{Definition}
\newtheorem{rem}[thm]{Remark}
\newtheorem{rem-ntn}[thm]{Remark and Notation}

\newenvironment{prf}{{\noindent \textbf{Proof:}\ }}{\hfill $\Box$\\ \smallskip}

\numberwithin{equation}{section}

\title[Solving Inverse Sturm-Liouville Problems via Cauchy Problems]{Solving Inverse Dirac-weighted Sturm-Liouville Problems via Cauchy problems}

\author[Min Zhao, Jiangang Qi and Xiao Chen]{ Min Zhao$^{a)}$, Jiangang Qi$^{b)}$ and Xiao Chen$^{c)}$ \\
Department of Mathematics, Shandong University, Weihai 264209, P.R. China \\
$a)$ \emph{Email}: \texttt{zhaomin215@mail.sdu.edu.cn}  \\
$b)$ \emph{Email}: \texttt{qjg816@163.com}  \\
$c)$ Author to whom correspondence should be addressed: \texttt{chenxiao@sdu.edu.cn}}

\begin{document}
\maketitle

\begin{abstract}
Building on the point interaction method developed in our previous work, this paper studies inverse eigenvalue problems for regular Sturm-Liouville problems with Dirac weights.
More precisely, we explicitly reconstruct the potential of regular Sturm-Liouville problems with the single-point and two-point Dirac weights from the solutions of a class of Cauchy problems which is completely determined by the first eigenvalues of a family of perturbed problems originating from moving point interaction models in quantum mechanics. Finally, the multi-point Dirac weighted case is also discussed.

\vspace{04pt}


\end{abstract}


\section{Introduction and problem statement}\label{sec:intro-prob-stat}
\medskip

In the Sturm-Liouville (abbreviate as S-L) theory, many eigenvalue problems often involve Dirac delta functions, although the coefficients of the classical Sturm-Liouville equations are considered in the space $L^1[a,b]$. More to the point, in physics, it is very important and common when the coefficients of second-order ordinary differential equations contain Dirac distributions (cf. \cite{AGHH1988, Car1995, CH1989, KP1931, Yang1967, McG1965, McG1966, Guli2019, CQ2021} etc.).

In \cite{CQ2021}, we have studied a forced vibration system of a string attached with finitely many particles, which can be described as a class of {\bf Dirac-weighted Sturm-Liouville eigenvalue problems}, in which the potentials are Lebesgue integrable functions and the weights are linear combinations of finitely many Dirac $\delta$-functions, namely, an ordinary differential equation with Dirichlet boundary conditions:
\begin{equation}\label{eqn:main-prob}
(P_n): \quad  -(pv^{\prime})^{\prime}(x)+q(x)v(x)=\lambda \sum_{j=1}^{n} m_j \delta (x -x_j) v(x),\ x\in[0,1],\ v(0)=0=v(1),
\end{equation}
where $0<x_0<\cdots<x_n<1$, $q\in L^1[0,1]$, $m_j\in \mathbb{R},\ j=1\cdots n,\ n\in \mathbb{N}$, $v(x)$ is the amplitude of the vibration, and $\delta(x-t)$ is the {\bf Dirac $\delta$-function} at $t\in(0,1)$ defined by
\begin{equation}\label{eqn:delta-funct}
\delta(x-t)=
\begin{cases}
\infty, &x=t,\\
          0, &x\neq t
\end{cases}
\ \text{and}\ \int_I\delta(x-t)\, dx=1,\ \ \forall\ I\subset[0,1]
\ \text{and}\ t\in I.
\end{equation}
Note that $$\int_0^1 f(x)\delta(x-t)\,{\rm d}x=\int_{t-\epsilon}^{t+\epsilon} f(x)\delta(x-t)\,{\rm d}x=f(t),$$ for any continuous function $f$ and $0<\epsilon<\min\{t,1-t\}$.

The solution of \eqref{eqn:main-prob} can be defined in the following way (cf.  \cite[(2.12)]{CQ2021}).
\begin{defn}\label{defn:solution-mainprob}
A function $y$ is called a {\bf solution of \eqref{eqn:main-prob}} if it satisfies
\begin{equation}
y\in\left\{z: z\in AC[0,1],\ z'\in\bigcup_{i=0}^{n} AC(I_j),\ \exists\ z'(x_j^{\pm}),\ j=1,\dots, n \right\},
\end{equation}
and
\begin{equation}
\left\{\aligned
&-y''(x)+q(x)y(x)=0,\ x\neq x_j,\ x\in[0,1],\\
&y(x_j^{-})=y(x_j^{+}),\ y'(x_j^{-})-y'(x_j^{+})=\lambda m_j y(x_j),\ j=1,\dots,n,\\
&y(0)=0=y(1).
\endaligned\right.
\end{equation}
where $I_j=(x_j, x_{j+1}),\  j=0,\dots, n,\ x_0=0,\ x_{n+1}=1$. $y'(x_j^{+}),\ y'(x_j^{-})$  respectively denote the classical  right-derivative and left-derivative, and $y'(0^{-}),\ y'(1^{+})$  is written as $y'(0),\ y'(1)$. $AC(I_j)$ denotes the space of absolutely continuous real-valued functions on an subinterval $I_j$.
\end{defn}

For the vibration system in \cite{CQ2021}, we consider a special external force determined by a non-homogeneous intensity $q(x)$ along $x$-axis and the perpendicular deflection $u(x,t)$, i.e., $F(x,t)=q(x)u(x,t)$. Hence, the integrable potential function $q(x)$ is viewed as the stationary potential of this non-uniform force field $F$. With the positions and masses of the particles on the string fixed, the potential $q(x)$ markedly alters the resulting vibrational frequencies. This naturally motivates a famous inverse spectral problem: \emph{if the vibrational frequencies of the string (i.e., the eigenvalues of the problem $(P_n)$), are known, can the potential $q(x)$ be uniquely recovered?}
This kind of inverse problem is of critical importance in structural dynamics and engineering diagnostics (cf. \cite{M1997, MD2002, HM2013, S2015} etc.).

On the one hand, in \cite{CQ2021}, we generalized Zhang's results \cite{ZWMQX2018} to the present case, and obtained complete results about the number of its Dirichlet eigenvalues of the problem $(P_n)$.
More precisely, we proved that, for any Sturm-Liouville equation with a non-zero integrable potential,  if the weight is a positive linear combination of $n$ Dirac $\delta$-functions, then the set of its Dirichlet eigenvalues either consists of at most $n$  (may be less than $n$, or even be $0$) different real numbers, or is the whole complex plane.
However, among the existing methods, for recovering the unknown potentials on the whole interval, we usually need infinitely many eigenvalues that are pre-given (cf. \cite{PT1987, FY2008, GW2014, FWW2016} etc.). Hence the finite Dirichlet eigenvalues of the problem $(P_n)$ are not enough to recover the potential uniquely.

On the other hand, in \cite{ZHQC2025, ZHQC2026}, we proposed a new approach to solving regular and singular inverse S-L problems based on  {\bf moving point interaction}, which originates from the {\bf $\delta$-point interaction} in quantum physics, the reader may refer to \cite{McG1965, McG1966, Yang1967, AGHH1988} and the references therein.
In the present paper, we will continue to adopt the strategy to tackle the inverse problem of $(P_n)$.

More precisely, for the S-L eigenvalue problem in \eqref{eqn:main-prob}, we produce a moving Dirac $\delta$-function perturbation similar to the procedures in \cite{ZHQC2025, ZHQC2026}. For fixed coupling constant $r>0$ and a variable perturbation position $t\in(0,1)$, we consider the following perturbed problem
\begin{equation}\label{eqn-interact-probl}
-y''(x)+[q(x)-r\delta(x-t)]y(x)=\lambda w(x)y(x), \ x\in(0,1),\ y(0)=0=y(1),
\end{equation}
where $q,\ w=\sum_{j=1}^{n} m_j \delta (x -x_j)$ are just the ones in \eqref{eqn:main-prob}.
As  $t$ runs over all real numbers in $(0,1)$, we obtain a family of {\bf perturbed problems} associated with the original problem \eqref{eqn:main-prob}, defined as \eqref{eqn-interact-probl} and denoted by ($P_n^t$).

\begin{rem}
If the perturbation is placed exactly at an endpoint, i.e. $t=0$ or $t=1$, the perturbation term $-r\delta(x-t)$ has no influence on the original problem ($P_0$).
\end{rem}

It is worth noticing that, compared with \cite{ZHQC2025}, the paper does not deal with S-L problems with integrable weights, but with Dirac $\delta$-function weights. And the coupling constant $r$ here is fixed, while the $r$ in \cite{ZHQC2025}  takes all the values in $[0, 1]$, so the required spectral data here are much less than  those in  \cite{ZHQC2025}.

Furthermore, later in this paper, we will solve our inverse problems via a series of Cauchy problems for first-order ordinal differential equations (abbreviate as ODEs), which are determined solely by the first eigenvalues of perturbed problems.

\begin{defn}
Suppose that every perturbed problem $(P_n^t)$ has a first eigenvalue. Denote  the first eigenvalue of $(P_n^t)$ for any $t\in(0,1)$ by $\lambda(t)$, which can be regarded as a function of the perturbation position $t$.
Then we call $\lambda(t)$ {\bf the first eigenvalue function} of the original problem $(P_n)$.
\end{defn}

Assume the problem
\begin{equation}\label{equ-w-1}
-y''+qy=\lambda y,\ y(0)=0=y(1)
\end{equation}
is left-definite,
where $q$ is just the one in \eqref{eqn:main-prob}.
Let $L(x)$ and $R(x)$ be the solutions of the equation $-y''+qy=0$ satisfying the initial conditions
\begin{equation}\label{ini-sols}
L(0)=0, \ L'(0)>0;\  R(1)=0, \ R'(1)<0.
\end{equation}
\begin{rem}\label{rem-left-defi}
The condition that \eqref{equ-w-1} is left-definite is equivalent to the statement that all eigenvalues of \eqref{equ-w-1} are positive (cf. \cite[Definition~5.2.1]{Zettl2005}), and the
following condition holds:
$$
{\bf (D_1)}: \quad  L(x)>0,\ R(x)>0,\ x\in (0,1).
$$
Since $0$ is not an eigenvalue of \eqref{equ-w-1}, we can assume that the Wronskian $W[L,R](x)=LR'-L'R\equiv -1$.

If \eqref{equ-w-1} is left-indefinite, all complex numbers may be eigenvalues of ($P_n$) (cf. \cite[Proposition~3.10]{CQ2021}). In this case, we replace ($P_n$) and ($P_n^t$) with the perturbed problems
\begin{equation}
\left\{\aligned
&-y''(x)+[q(x)+c]y(x)=\lambda wy(x),\  x\in[0,1],\\
&  y(0)=0=y(1),
\endaligned\right.
\end{equation}
and
\begin{equation}
\left\{\aligned
&-y''(x)+[q(x)+c-r\delta(x-t)]y(x)=\lambda wy(x),\  x\in[0,1],\\
&  y(0)=0=y(1),
\endaligned\right.
\end{equation}
respectively, where $q,\ w=\sum_{j=1}^{n} m_j \delta (x -x_j)$ are just the ones in \eqref{eqn:main-prob},  $c>\widehat{\lambda}_1$ is a constant, and $\widehat{\lambda}_1$ is the first eigenvalue of the problem \eqref{equ-w-1}.
\end{rem}

To characterize $\lambda(t)$ in the next section, we introduce two auxiliary solutions $\varphi(x,\lambda)$ and $\psi(x,\lambda)$ of the equation
$$
-y''+qy=\lambda w y,
$$
in \eqref{eqn:main-prob} satisfying the initial conditions
\begin{equation}\label{ini-sols-lamb}
 \varphi(0)=0, \ \varphi'(0)=L'(0);\  \psi(1)=0, \ \psi'(1)=R'(1).
\end{equation}

Analogous to \cite{ZHQC2025, ZHQC2026}, the following discriminant condition of an eigenvalue of \eqref{eqn:main-prob} is also necessary.
\begin{thm}(cf. \cite[Lemma 3.3]{ZHQC2025})\label{thm-equi}
$\lambda$ is an eigenvalue of \eqref{eqn:main-prob} if and only if it satisfies
\begin{equation}\label{eqn:equiv-1st-eigenv}
 r\varphi(t,\lambda)\psi(t,\lambda)=-W[\varphi,\psi],\ r>0,\ t\in(0,1),
\end{equation}
where $\varphi, \psi$ are defined as in \eqref{ini-sols-lamb}.
\end{thm}

The paper is organized as follows. Section~\ref{sec:single-dirac} deals with the case of single-point Dirac weight, providing the complete reconstruction method. Section~\ref{sec:two-dirac} extends this analysis to the case of two-point Dirac delta weight.
In the final section, we also discuss the general multi-point case, and propose several open problems for future research.

\bigskip

\section{For the case of single-point Dirac weight}\label{sec:single-dirac}
\medskip

In this section, we consider the S-L problem \eqref{eqn:main-prob} with a single Dirac delta weight, i.e. the weight has the form $w(x)=w_0\delta(x-a)$, where $w_0>0,\ a\in(0,1)$. The original problem, denoted by ($P_1$), is
\begin{equation}\label{sing-dirac-weig}
\left\{\aligned
&-y''(x)+q(x)y(x)=\lambda w_0\delta(x-a)y(x),\  x\in[0,1],\\
&  y(0)=0=y(1).
\endaligned\right.
\end{equation}

The corresponding perturbed problem ($P_1^t$)
\begin{equation}\label{per-sing-dirac-weig}
\left\{\aligned
&-y''(x)+[q(x)-r\delta(x-t)]y(x)=\lambda w_0\delta(x-a)y(x),\  x\in[0,1],\\
&  y(0)=0=y(1),
\endaligned\right.
\end{equation}
where $r>0$ is fixed, and $t\in [0,1]$ is the movable perturbation position.  Clearly, $(P_1^1)=(P_1^0)=(P_1).$

\begin{rem}\label{rem-sing-ori-exis}
From Remark \ref{rem-left-defi}, the condition that \eqref{equ-w-1} is left-definite is equivalent to ${\bf (D_1)}$ holds.
Under the condition ${\bf (D_1)}$, the hypotheses $(H_0)$ and $(H)$ in \cite{CQ2021} hold, and it follows from \cite[Proposition 3.7]{CQ2021} that the problem ($P_1$) has a unique positive eigenvalue, denoted by $\lambda_1$.
\end{rem}

The following lemma will play a central role for the reconstruction result.
\begin{lem}\label{lem:basic-relation}
Suppose that \eqref{equ-w-1} is left-indefinite, then
\begin{equation}\label{sing-weig-first-eig}
w_0 L(a)R(a)\lambda(0)=1,
\end{equation}
where  $L,\ R$ are defined in \eqref{ini-sols}.
\end{lem}

\prf From Remark \ref{rem-sing-ori-exis}, we have $\lambda(0)=\lambda_1$.
By Definition \ref{defn:solution-mainprob}, the solutions $\varphi(x,\lambda)$ and $\psi(x,\lambda)$ introduced in \eqref{ini-sols-lamb} can be chosen as
\begin{equation}\label{equ-secle-1-var}
\varphi(x,\lambda)=\left\{\aligned
&L(x),\ &x\in[0,a],\\
&\alpha_1 L(x)+\beta_1 R(x),\  &x\in[a,1],
\endaligned\right.
\end{equation}
\begin{equation}\label{equ-secle-1-psi}
\psi(x,\lambda)=\left\{\aligned
&\mu_1 L(x)+\gamma_1 R(x),\  &x\in[0,a],\\
&R(x),\ &x\in[a,1],
\endaligned\right.
\end{equation}
where $\alpha_1=\alpha_1(\lambda)$, $\beta_1=\beta_1(\lambda)$, $ \mu_1=\mu_1(\lambda)$ and $\gamma_1=\gamma_1(\lambda)$.

From the continuity of $\varphi$ and the jump condition  for its derivative at $x=a$, we have
$$
L(a)=\alpha_1 L(a)+\beta_1 R(a),  \ L'(a)-[\alpha_1 L'(a)+\beta_1 R'(a)]=\lambda w_0 L(a).
$$
Using $W[L,R](x)\equiv -1$, we obtain
\begin{equation}\label{alpha1-beta1}
 \alpha_1(\lambda)=1-w_0 L(a)R(a)\lambda, \ \beta_1(\lambda)=w_0L^2(a)\lambda.
\end{equation}
Similarly, from the continuity of $\psi$ and the jump condition  for its derivative at $x=a$, we also have
\begin{equation}\label{mu1-gamma1}
 \mu_1(\lambda)=w_0R^2(a)\lambda, \ \gamma_1(\lambda)=1-w_0L(a)R(a)\lambda=\alpha_1(\lambda).
\end{equation}
Since $\lambda(0)=\lambda_1$ is the first eigenvalue of the problem ($P_1$), one can see that $\varphi(1,\lambda(0))=0$. From the expression
$$
\varphi(1,\lambda(0))= \alpha_1(\lambda(0))L(1)+\beta_1(\lambda(0))R(1)
$$
and the fact that $R(1)=0$, we have $ \alpha_1(\lambda(0))=0$. Then, the equation \eqref{sing-weig-first-eig} follows from \eqref{alpha1-beta1}.
\qed

\begin{rem} We have $\lambda(a)=0$ if and only if $r=\lambda_1 w_0$. For the case $r\neq\lambda_1 w_0$ and
special perturbation position $t=a$, one has that $\lambda(a)$ is the first eigenvalue of the problem
$$
-y''+[q-r\delta(x-a)]y=\lambda w_0\delta(x-a) y,\ y(0)=0=y(1),
$$
or equivalently, $\lambda(a)$ is the first eigenvalue of
$$
-y''+qy=\lambda \left(w_0+\frac{r}{\lambda(a)}\right)\delta(x-a)y,\ y(0)=0=y(1).
$$
Applying Lemma~\ref{lem:basic-relation} yields
$$
\left(w_0+\frac{r}{\lambda(a)}\right) L(a)R(a)\lambda(a)=1,
$$
which, combined with equation \eqref{sing-weig-first-eig}, implies that $\lambda(a)=\lambda_1-\frac{r}{w_0}$.
\end{rem}

\begin{lem}\label{lem-exist-prop-lamd1-sing}
Suppose that \eqref{equ-w-1} is left-indefinite and $$ 0<r<\min\left\{ \min\limits_{t\in(0,a]}\frac{1}{(R-\frac{R(a)}{L(a)}L)L}, \min\limits_{t\in[a,1)}\frac{1}{(L-\frac{L(a)}{R(a)}R)R} \right\}.$$
Then the first eigenvalue function  $\lambda(t)$ of ($P_1$) exists and satisfies\\
$(i)$ $\lambda''(t)$ exists a.e. on $(0,1)$;\\
$(ii)$ $\lambda(t)<\lambda(0)=\lambda(1)=\lambda_1,\ t\in (0,1)$.
\end{lem}
\prf For any $t\in(0,1)$, by Theorem \ref{thm-equi}, we know that $\lambda(t)$ is an eigenvalue of $(P_1^t)$ if and only if it satisfies
\begin{equation}\label{sing-weig-equiv}
r\varphi(t,\lambda(t))\psi(t,\lambda(t))=-W[\varphi,\psi],
\end{equation}
where $\varphi$ and $\psi$ are chosen as \eqref{equ-secle-1-var} and \eqref{equ-secle-1-psi}, respectively.
We have
$$
W[\varphi,\psi]=\left\{\aligned
& W[L, \mu_1 L+\gamma_1 R]=\gamma_1 W[L,R]=-\gamma_1,\  &x\in[0,a],\\
& W[\alpha_1 L+\beta_1 R, R]=\alpha_1W[L,R]=-\alpha_1, \ &x\in[a,1],
\endaligned\right.
$$
where $\alpha_1=\alpha_1(\lambda(t))$, $\beta_1=\beta_1(\lambda(t))$, $\mu_1=\mu_1(\lambda(t))$ and $\gamma_1=\gamma_1(\lambda(t))$.

For $t\in[a,1)$, the equation \eqref{sing-weig-equiv} becomes
\begin{equation}\label{eq:tina1}
r[\alpha_1 L(t)+\beta_1 R(t)]R(t)=\alpha_1.
\end{equation}
Since $W[L,R](x)\equiv -1$, it holds that
\begin{equation}\label{eq:R-representation}
R(t)=L(t)\int^1_t\frac{{\rm d}s}{L^2(s)}.
\end{equation}
Substituting \eqref{eq:R-representation} into \eqref{eq:tina1} gives
\begin{equation}\label{equ-sing-LR}
r\left[\alpha_1+\beta_1\int^1_t\frac{{\rm d}s}{L^2(s)}\right]L^2(t)\int^1_t\frac{{\rm d}s}{L^2(s)}=\alpha_1.
\end{equation}
Define
\begin{equation}\label{equ-sing-U-def-U}
U(t)=\frac{L(a)}{R(a)}\int^1_t\frac{{\rm d}s}{L^2(s)},
\end{equation}
then we have
$$U(a)=\frac{L(a)}{R(a)}\int^1_a\frac{{\rm d}s}{L^2(s)}=\frac{R(a)}{R(a)}=1,$$
and from \eqref{equ-sing-LR}, we deduce that $U$ satisfies the differential equation
\begin{equation}\label{equ-sing-differ-equt-U}
\alpha_1U'(t)=-rU\left[\alpha_1+\beta_1\frac{R(a)}{L(a)}U\right].
\end{equation}
Note that \eqref{sing-weig-first-eig}, together with \eqref{alpha1-beta1}, implies
\begin{equation}\label{equ-single-beta1-alph1}
\beta_1\frac{R(a)}{L(a)}=L(a)R(a)\lambda=\frac{\lambda(t)}{w_0\lambda(0)},\
\alpha_1=1-\lambda(t)w_0L(a)R(a)=1-\frac{\lambda(t)}{\lambda(0)}.
\end{equation}
Therefore, the perturbed problem $(P_1^t)$ has an eigenvalue if and only if, for fixed $t$, the equation
\begin{equation}\label{equ-sing-lamt}
\lambda(t)(rU(1-U)+U')=\lambda(0)(U'+rU)
\end{equation}
has a solution $\lambda(t)$.

We first prove the existence of the eigenvalue of $(P_1^t)$. Since $r< \min\limits_{t\in[a,1)}\frac{1}{(L-\frac{L(a)}{R(a)}R)R} $, we have $[rU(1 - U) + U'](t)< 0$ for all $t\in[a,1)$ and the perturbed problem $(P_1^t)$ has an eigenvalue $\lambda(t)$ for any $t\in[a,1)$.

Next, we prove $\lambda(t)<\lambda(0),\ t\in[a,1)$. We consider three cases as follows:

\noindent Case $1$. If $[U' + rU](t)=0$, then $\lambda(t)=0<\lambda(0)$.

\noindent Case $2$. If $[U' + rU](t)>0$, then $\lambda(t)(rU(1-U)+U')>0$. Since $[rU(1-U)+U'](t)<0$, it follows that $\lambda(t)<0<\lambda(0)$.

\noindent Case $3$. If $[U' + rU](t)<0$, then $\lambda(t)>0$.
By \eqref{equ-sing-lamt}, we have
$$\lambda(t)=\lambda(0)\left(1+\frac{rU^2}{U'+rU-rU^2}\right)<\lambda(0).$$
Finally, by the definition \eqref{equ-sing-U-def-U} of $U$ and the equation \eqref{equ-sing-lamt}, we can see that $\lambda''(t)$ exists a.e. on $(a,1).$

For $t\in(0,a]$, the condition \eqref{sing-weig-equiv} takes the form
\begin{equation}\label{eq:region1-main}
rL(t)[\mu_1 L(t)+\gamma_1 R(t)]=\gamma_1.
\end{equation}
Using
$$
L(t)=R(t)\int^t_0\frac{{\rm d}s}{R^2(s)},
$$
which is a consequence of $W[L,R] \equiv -1$,
the equation \eqref{eq:region1-main} becomes
\begin{equation}\label{equ-sing-V-only}
r\left[\mu_1\int^t_0\frac{{\rm d}s}{R^2(s)}+\gamma_1\right]R^2(t)\int^t_0\frac{{\rm d}s}{R^2(s)}=\gamma_1.
\end{equation}
Set
\begin{equation}\label{equ-sing-def-V}
V(t)=\frac{R(a)}{L(a)}\int^t_0\frac{{\rm d}s}{R^2(s)},
\end{equation}
then
$$
V(a)=\frac{R(a)}{L(a)}\int^a_0\frac{{\rm d}s}{R^2(s)}=\frac{L(a)}{L(a)}=1,
$$
and from \eqref{equ-sing-V-only}, we deduce that $V$ satisfies the differential equation
\begin{equation}\label{equ-sing-diffe-equ-V}
\gamma_1V'(t)=rV\left[\gamma_1+\mu_1\frac{L(a)}{R(a)}V\right].
\end{equation}
Note that \eqref{mu1-gamma1}, together with \eqref{sing-weig-first-eig},  implies
\begin{equation}\label{equ-sing-mu1-gamma1}
\mu_1\frac{L(a)}{R(a)}=\beta_1\frac{R(a)}{L(a)}=\frac{\lambda(t)}{w_0\lambda(0)},\
\gamma_1=\alpha_1=1-\frac{\lambda(t)}{\lambda(0)}.
\end{equation}
Therefore, the perturbed problem $(P_1^t)$ has an eigenvalue if and only if, for fixed $t$, the equation
\begin{equation}\label{equ-sing-lamt-0-a}
\lambda(t)(rV(V-1)+V')=\lambda(0)(V'-rV).
\end{equation}
has a solution $\lambda(t)$.
Since $r< \min\limits_{t\in(0,a]}\frac{1}{(R-\frac{R(a)}{L(a)}L)L}$, we have $[rV(V-1)+V'](t)>0$ for all $t\in(0,a]$.
Similar to the case that $t\in [a,1)$, we can also prove that
$\lambda(t)$ exists and $\lambda(t)<\lambda(0)$.
Moreover, by the definition \eqref{equ-sing-def-V} of $V$ and the equation \eqref{equ-sing-lamt-0-a}, we have that $\lambda''(t)$ exists a.e. on $(0,a).$


\qed

The following theorem states the uniqueness and reconstruction result for $q$ in $(P_1)$.
Based on the first eigenvalue of the perturbed problem ($P_1^t$),
we first derive the expressions \eqref{sing-weig-right-sol} and \eqref{sing-weig-left-sol} for the solutions $L$ and $R$ defined in \eqref{ini-sols}, respectively. These expressions then enable the reconstruction of the potential.

\begin{thm}\label{thm-sing-recov} Under the conditions in Lemma \ref{lem-exist-prop-lamd1-sing}, let $\lambda(t)$ be the first eigenvalue function of ($P_1$). Then $q$ can be uniquely determined by $\lambda(t)$, and
$$q(x)=\left\{\aligned
&R''(x)/R(x),\ x\in(0,a],\\
&L''(x)/L(x),\ x\in[a,1),
\endaligned\right.$$
where $L$ and $R$ are defined as in \eqref{ini-sols} and can be expressed explicitly by
\begin{equation}\label{sing-weig-left-sol}
R(t)=\left\{\frac{r\left[e^{r(a-t)}+r \int^{a}_{t} \frac{\lambda(s)}{\lambda_{1}-\lambda(s)} e^{r(s-t)} \mathrm{d}s+\frac{\lambda(t)}{\lambda_{1}-\lambda(t)}\right]}
{\left[e^{r(a-t)}+r \int^{a}_{t} \frac{\lambda(s)}{\lambda_{1}-\lambda(s)} e^{r(s-t)} \mathrm{d}s\right]^{2}}\right\}^{-1/2},\ t\in(0,a],
\end{equation}
\begin{equation}\label{sing-weig-right-sol}
L(t)=\left\{
\frac{r\left[e^{r(t-a)}+r \int^{t}_a \frac{\lambda(s)}{\lambda_{1}-\lambda(s)} e^{r(t-s)} \mathrm{d}s+\frac{\lambda(t)}{\lambda_{1}-\lambda(t)}\right]}
{\left[e^{r(t-a)}+r \int^{t}_a \frac{\lambda(s)}{\lambda_{1}-\lambda(s)} e^{r(t-s)} {\rm d}s\right]^{2}}
\right\}^{-1/2},\ t\in[a,1).
\end{equation}
\end{thm}

\prf For $t\in[a,1)$, it follows from \eqref{equ-sing-differ-equt-U} and \eqref{equ-single-beta1-alph1} that $U$ satisfies the Cauchy problem
\begin{equation}\label{equ-single-cauch-U}
\left\{\aligned
& U'(t)=-rU\left[1+\frac{\lambda(t)}{\lambda(0)-\lambda(t)}U\right],\ t\in[a,1),\\
& U(a)=1.
\endaligned\right.
\end{equation}
Since all the coefficients in \eqref{equ-single-cauch-U} are uniquely determined by $\lambda(t)$,  Lemma \ref{lem-exist-prop-lamd1-sing}, together with the standard theory of ODEs, guarantees the existence and uniqueness of the solution $U$.
Set $z=U^{-1}$, then $z$ satisfies the initial value problem
\begin{equation}
\left\{\aligned
& z'(t)=rz+\frac{r\lambda(t)}{\lambda(0)-\lambda(t)},\ t\in[a,1),\\
& z(a)=1,
\endaligned\right.
\end{equation}
whose solution is
\begin{equation}\label{eq:z-solution}
z(t)=e^{r(t-a)}+r \int^{t}_a \frac{\lambda(s)}{\lambda_{1}-\lambda(s)} e^{r(t-s)} \mathrm{d}s.
\end{equation}
Then, from the definition of $U$ in \eqref{equ-sing-U-def-U}, we have
\begin{equation}\label{equ-L-U-z}
L(t)=\left(L(a)/R(a)\right)^{1/2}\left(-U'(t)\right)^{-1/2}
=\left(L(a)/R(a)\right)^{1/2}\left(z'/z\right)^{-1/2}(t),
\end{equation}
and substituting \eqref{eq:z-solution} into \eqref{equ-L-U-z} gives the reconstruction formula \eqref{sing-weig-right-sol}.
Finally, by the definition of $L$ in \eqref{ini-sols}, we obtain
\begin{equation}\label{equ-sing-a-1-qL}
q(x)=L''(x)/L(x),\ x\in[a,1).
\end{equation}

For $t\in(0,a]$, it follows from \eqref{equ-sing-diffe-equ-V} and \eqref{equ-sing-mu1-gamma1} that $V$ satisfies the Cauchy problem
\begin{equation}\label{equ-single-cauch-V}
\left\{\aligned
& V'(t)=rV\left[1+\frac{\lambda(t)}{\lambda(0)-\lambda(t)}V\right],\ t\in(0,a],\\
& V(a)=1.
\endaligned\right.
\end{equation}
Using Lemma \ref{lem-exist-prop-lamd1-sing} again, we also can obtain the existence and uniqueness of the solution $V$, which is uniquely determined by $\lambda(t)$.
Set $y=V^{-1}$, then $y$ satisfies
\begin{equation}
\left\{\aligned
& y'(t)=-ry-\frac{r\lambda(t)}{\lambda(0)-\lambda(t)},\ t\in(0,a],\\
& y(a)=1,
\endaligned\right.
\end{equation}
whose solution is
\begin{equation}\label{eq:y-solution}
y(t)=e^{r(a-t)}+r \int^{a}_x \frac{\lambda(s)}{\lambda_{1}-\lambda(s)} e^{r(s-t)} \mathrm{d}s.
\end{equation}
Then from the definition of $V$ in \eqref{equ-sing-def-V}, we have
\begin{equation}\label{equ-sing-R-V-y}
R(t)=\left(R(a)/L(a)\right)^{1/2}\left(V'(t)\right)^{-1/2}
=\left(R(a)/L(a)\right)^{1/2}\left(-y'/y\right)^{-1/2}(t),
\end{equation}
and substituting \eqref{eq:y-solution} into \eqref{equ-sing-R-V-y} yields exactly the reconstruction formula \eqref{sing-weig-left-sol}.
Finally, by the definition of $R$ in \eqref{ini-sols}, we obtain
\begin{equation}\label{equ-sing-0-a-qR}
q(x)=R''(x)/R(x),\ x\in(0,a].
\end{equation}
 \qed
\begin{rem}\label{rem-sing-dira}
Based on the relations between $U,\ V$ and $L,\ R$  given in \eqref{equ-L-U-z} and \eqref{equ-sing-R-V-y}, together with \eqref{equ-sing-a-1-qL} and \eqref{equ-sing-0-a-qR}, the potential $q$ can also be expressed as  $$q(x)=\left\{\sqrt{1/Y'(x)}\right\}''\Big/\!\sqrt{1/Y'(x)}, \ x\in(0,1),$$
where $Y$ is defined as $Y(x)=\left\{\aligned
& V(x),\ x\in(0,a],\\
& U(x),\ x\in[a,1),
\endaligned\right.$
and
$Y$ satisfies the following {\bf Cauchy problem}:
\begin{equation}\label{3212}
\left\{\aligned
& Y'(t)=ra_i\left(1+b_i Y\right),\ t\in [x_{i-1},x_i],\\
& Y(t_i)=1,
\endaligned\right.
\end{equation}
where $i=1, 2$, $t_1=t_2=a$, $x_0=0,\ x_1=a,\ x_2=1$, and the coefficients
$$a_1=-1,\ a_2=1,\ b_1=b_2=\frac{\lambda(t)}{\lambda(0)-\lambda(t)}$$
are given by \eqref{equ-single-cauch-V} and \eqref{equ-single-cauch-U}.
\end{rem}

\begin{rem}\label{rem-comp}
Theorem \ref{thm-sing-recov} recovers $q$ using $\lambda(t)$ that depends only on the adjusted perturbed position $t$.
By comparison, \cite[Theorem~4.1]{ZHQC2025} reconstructs $q$ using the first eigenvalue function $\lambda(t,r)$ of the problem
\begin{equation}\label{pro-w-integ}
-y''(x)+q(x)y(x)=\lambda w(x)y(x), \ x\in(0,1),\ y(0)=0=y(1),
\end{equation}
with $0<w\in L^1[0,1]$, by varying the interaction position $t$ along the open interval $(0,1)$ as well as adjusting the intensity $r$ of the perturbed problem
$$-y''(x)+[q(x)-r\delta(x-t)]y(x)=\lambda w(x)y(x), \ x\in(0,1),\ y(0)=0=y(1).$$
The reconstruction formula is
\begin{equation}\label{eqn:reconstruct-formula}
q(x)=\frac{\varphi''_0(x)}{\varphi_0(x)}+\lambda_1 w(x),\
\varphi_0(x)=\sqrt{-\frac{\partial \lambda(x,0)}{\partial r}},\ x\in(0,1),
\end{equation}
where $\lambda_1=\lambda(t,0)$, is the first eigenvalue of \eqref{pro-w-integ}.

If the positive parameter $r$ in ($P_1^t$) is also allowed to vary, then we denote by $\lambda(t,r)$ the first eigenvalue of $(P_1^t)$ for any $t\in(0,1)$ and $r>0$, and $\lambda(t,r)$ can be viewed as a function of $t$ and $r$. Note that $w(x)=w_0\delta(x-a)$ here, so $w\equiv 0$ on both $[0,a)$ and $[a,1)$.

For $t\in(a,1)$, it follows from \eqref{equ-sing-lamt} that
$$\frac{\partial \lambda(t,0)}{\partial r}=-\lambda(0)\frac{L(a)}{R(a)}R^2.$$
Substituting this identity into \eqref{eqn:reconstruct-formula} yields  $q=R''/R$, which is exactly consistent with the conclusion stated in Theorem \ref{thm-sing-recov}.

For $t\in(0,a)$, by virtue of  \eqref{equ-sing-lamt-0-a}, we obtain
$$\frac{\partial \lambda(t,0)}{\partial r}=-\lambda(0)\frac{R(a)}{L(a)}L^2,$$
and plugging this into \eqref{eqn:reconstruct-formula} further yields $q=L''/L$, which also agrees with Theorem \ref{thm-sing-recov}.

Consequently, Theorem \ref{thm-sing-recov} can be viewed as a generalization of  \cite[Theorem~4.1]{ZHQC2025} to the setting with single-point Dirac weights. Moreover, it provides a more explicit formulation of $q$, and requires fewer spectral data.
\end{rem}

\bigskip

\section{For the case of two-point Dirac weight}\label{sec:two-dirac}
\medskip

In this section, we  extend the reconstruction method to the case where the weight consists of two-point Dirac delta distribution  $$w(x)=m_1\delta(x-x_1)+m_2\delta(x-x_2),\ m_1,m_2>0,\ 0=x_0<x_1<x_2<x_3=1.$$
The corresponding unperturbed eigenvalue problem, denoted by $(P_2)$, is
\begin{equation}\label{two-dirac-weig}
\left\{\aligned
&-y''(x)+q(x)y(x)=\lambda [m_1\delta(x-x_1)+m_2\delta(x-x_2)]y(x),\  x\in[0,1],\\
& y(0)=0=y(1),
\endaligned\right.
\end{equation}
and the corresponding perturbed problem, denoted by ($P^t_2$), is
\begin{equation}\label{two-dirac-weig-per}
\left\{\aligned
&-y''(x)+[q(x)-r\delta(x-t)]y(x)=\lambda [m_1\delta(x-x_1)+m_2\delta(x-x_2)]y(x),\  x\in[0,1],\\
& y(0)=0=y(1),
\endaligned\right.
\end{equation}

\begin{rem}\label{rem-two-ori-exis}
Under the condition that \eqref{equ-w-1} is left-definite, we have  $$
{\bf (D_2)}: L_2R_1-L_1R_2\neq 0,
 $$
where $L_j=L(x_j),\ R_j=R(x_j),\ j=1,2.$ Otherwise, we can prove that $W[L,R]=0$, which contradicts that $W[L,R](x)\equiv -1$ in Remark \ref{rem-left-defi}.
Under the condition ${\bf (D_2)}$, problem ($P_2$) has exactly two eigenvalues (cf. \cite[Proposition 3.7]{CQ2021}).
\end{rem}
In this section, unless otherwise specified, we continue to use the notation established in Section 2, that is $\lambda_1$ and $\lambda(t)$ are the first eigenvalue and the first eigenvalue function of $(P_2)$,  respectively. And then  $\lambda(0)=\lambda(1)=\lambda_1$ holds.

Set $x_0=0,\ x_3=1$. Because the two Dirac masses divide the interval into three subintervals $I_j=(x_{j-1},x_j),\ 1 \leq j\leq 3$, the solutions $\varphi$ and $\psi$ defined in \eqref{ini-sols-lamb} can be written as
\begin{equation}\label{two-dirac-experss-varphi}
\varphi(x,\lambda)=\left\{\aligned
&L(x),\ &&x\in[0,x_1],\\
&\alpha_2 L(x)+\beta_2 R(x),\  &&x\in[x_1,x_2],\\
&\alpha_3 L(x)+\beta_3 R(x),\  &&x\in[x_2,1],
\endaligned\right.
\end{equation}
and
\begin{equation}\label{two-dirac-experss-psi}
\psi(x,\lambda)=\left\{\aligned
&\mu_1 L(x)+\gamma_1 R(x), &&x\in[0,x_1],\\
&\mu_2 L(x)+\gamma_2 R(x),\  &&x\in[x_1,x_2],\\
&R(x),\ &&x\in[x_2,1],
\endaligned\right.
\end{equation}
where $\alpha_i=\alpha_i(\lambda),\ \beta_i=\beta_i(\lambda)$, $2\leq i\leq 3$; $\mu_j=\mu_j(\lambda)$, $\gamma_j=\gamma_j(\lambda)$, $1\leq j\leq 2$.
For convenience, we set
$$
\alpha_1=1, \ \beta_1=0;\ \mu_3=0,\ \gamma_3=1,
$$
so that both of \eqref{two-dirac-experss-varphi} and \eqref{two-dirac-experss-psi} can be uniformly expressed as
$$\varphi(x,\lambda)=\alpha_j L(x)+\beta_j R(x),\  x\in I_j,$$
and
$$\psi(x,\lambda)=\mu_j L(x)+\gamma_j R(x),\  x\in I_j,$$ where $1\le j\le 3$.
Applying the continuity of $\varphi$ and the jump conditions for its derivative at $x_j$ for $j=1,2$, that is,
$$
\left\{\aligned
& \varphi(x_j-0)=\varphi(x_j+0),\\
& \varphi'(x_j-0,\lambda)-\varphi'(x_j+0,\lambda)=\lambda m_j\varphi(x_j,\lambda),
\endaligned\right.
$$
we obtain
\begin{equation}\label{eqn:alpha-beta-relation-two-pt-case}
\left\{\aligned
& \alpha_{j+1}=\alpha_j-\lambda m_j(\alpha_j L_j+\beta_j R_j)R_j, \ j=1,2,\\
& \beta_{j+1}=\beta_j+\lambda m_jL_j(\alpha_j L_j+\beta_j R_j),  \ j=1,2.
\endaligned\right.
\end{equation}
Substituting $\alpha_1=1$ and $\beta_1=0$ into \eqref{eqn:alpha-beta-relation-two-pt-case}, we obtain
\begin{equation}\label{euq-two-calu-alph-beta}
\left\{\aligned
&\alpha_2=1-\lambda m_1L_1R_1,\ \alpha_3=1-\lambda m_1L_1R_1-\lambda m_2R_2(\alpha_2L_2+\beta_2 R_2),\\
&\beta_2=\lambda m_1 L^2_1,\ \beta_3=\lambda m_1 L^2_1+\lambda m_2L_2(\alpha_2L_2+\beta_2 R_2).
\endaligned\right.
\end{equation}
By \eqref{euq-two-calu-alph-beta}, we have
\begin{equation}\label{equ-alph3}
\alpha_3(\lambda)=1-\lambda[m_1L_1R_1+m_2L_2R_2]+\lambda^2m_1m_2L_1R_2[L_2R_1-L_1R_2].
\end{equation}
Applying the continuity and the jump conditions of $\psi$ at $x_j$ for $j=1,2$, that is,
$$
\left\{\aligned
& \psi(x_j-0)=\psi(x_j+0),\\
& \psi'(x_j-0,\lambda)-\psi'(x_j+0,\lambda)=\lambda m_j\psi(x_j,\lambda),
\endaligned\right.
$$
we obtain
\begin{equation}\label{euq-two-calu-mu-gamma}
\left\{\aligned
&\mu_3=0,\ \mu_2= \lambda m_2R_2^2,\ \mu_1=\lambda(m_1R_1^2+m_2R_2^2)+\lambda^2 m_1m_2R_1R_2(L_1R_2-L_2R_1),\\
&\gamma_3=1,\ \gamma_2=1-\lambda m_2L_2R_2,\ \gamma_1=\alpha_3,
\endaligned\right.
\end{equation}

The following lemma establishes the existence and differentiability of the eigenvalue of perturbed problem ($P^t_2$).
\begin{lem}\label{lem-exis-prop-lamd1-two}
Suppose \eqref{equ-w-1} is left-definite. Let
\begin{equation}\label{equ-condi-exis}
0<r\neq\left\{\frac{W'(m_1X_1+m_2X_2+m_1m_2X_3)}{|c|}:c(t)\neq0,\ t\in(x_1,x_2]\right\}
\end{equation}
where
$$c(t)=m_1m_2[ -(1-W)^2X_1X_2-X_3W(1-W) ]-\left( m_1X_1+W^2m_2\frac{X_1X_2-X_3}{X_1}  \right)+W(m_1X_1+m_2X_2),$$ $W=\frac{R_1L}{L_1R},\ W'=\frac{R_1}{L_1R^2}$ and
\begin{equation}\label{two-var-x1-x2-x3}
X_1=L_1R_1,\ X_2=L_2R_2,\ X_3=L_1R_2[L_2R_1-L_1R_2].
\end{equation}
Then the first eigenvalue function  $\lambda(t)$ of ($P_2$) exists and satisfies $\lambda''(t)$ exists a.e. on $(0,1)$.
\end{lem}
\prf For any $t\in(0,1)$, by Theorem \ref{thm-equi}, we know that, $\widetilde{\lambda}(t)$ is an eigenvalue of ($P^t_2$) if and only if it satisfies
\begin{equation}\label{two-weig-equiv}
r\varphi(t,\widetilde{\lambda}(t))\psi(t,\widetilde{\lambda}(t))=-W[\varphi,\psi].
\end{equation}
where $\varphi$ and $\psi$ are chosen as in \eqref{two-dirac-experss-varphi} and \eqref{two-dirac-experss-psi}, respectively.
And we have
$$
W[\varphi,\psi]\equiv-\gamma_1=-\alpha_3=-(\alpha_2\gamma_2-\beta_2\mu_2),
$$
where $\alpha_j=\alpha_j(\widetilde{\lambda}(t)),\ \beta_j=\beta_j(\widetilde{\lambda}(t))$,  $\mu_j=\mu_j(\widetilde{\lambda}(t))$, $\gamma_j=\gamma_j(\widetilde{\lambda}(t))$, $1\leq j\leq 3$.

For $t\in(0,x_1]$, equation \eqref{two-weig-equiv} becomes
\begin{equation}\label{two-weig-0-x1-w}
r[\mu_1 L(t)+\gamma_1 R(t)]L(t)=\gamma_1.
\end{equation}
Using the relation
\begin{equation}\label{LR}
L(t)=R(t)\int^t_0\frac{{\rm d}s}{R^2(s)},\ t\in(0,1),
\end{equation}
\eqref{two-weig-0-x1-w} can be written as
\begin{equation}\label{two-weig-0-x1-LR-w}
r\left[\mu_1\int^t_0\frac{{\rm d}s}{R^2(s)}+\gamma_1\right]R^2(t)\int^t_0\frac{{\rm d}s}{R^2(s)}=\gamma_1.
\end{equation}
Define
\begin{equation}\label{two-weig-0-x1-defi-V}
V(t)=\frac{R_1}{L_1}\int^t_0\frac{{\rm d}s}{R^2(s)},
\end{equation}
then equation \eqref{two-weig-0-x1-LR-w} reduces to
\begin{equation}\label{two-weig-0-x1-equi-cont}
A_1 \lambda^2+B_1 \lambda+C_1=0,
\end{equation}
where $$A_1=m_1m_2X_3(V'-rV+rV^2),$$
$$ B_1=-(V'-rV)(m_1X_1+m_2X_2)-rV^2\left(m_1X_1+m_2\frac{X_1X_2-X_3}{X_1}\right),\ C_1=V'-rV.$$

We now prove the existence of the first eigenvalue of $(P_2^t)$. We divide three cases.

\noindent Case $1$. If $A_1=0$, i.e., $V'-rV+rV^2=0$, then $B_1=rV^2m_2\frac{X_3}{X_1}>0$, and hence $(P_2^t)$ has an eigenvalue.

\noindent Case $2$. If $A_1<0$, then $C_1<0$. From the fact that $X_1X_2>X_3>0$, it follows that the discriminant
\begin{align*}
\Delta_1=B_1^2-4A_1C_1&=\left[rV^2m_2\frac{X_3}{X_1}-(rV-V'-rV^2)(m_1X_1-m_2X_2)\right]^2\\
 &+4m_1m_2(rV-V'-rV^2)^2(X_1X_2-X_3)>0.
\end{align*}
Thus $(P_2^t)$ has two eigenvalues.

\noindent Case $3$. $A_1>0$.

\noindent If $C_1<0$, then $(P_2^t)$ has two eigenvalues from the property of quadratic function.

\noindent If $C_1>0$, using $X_1X_2>X_3>0$, we obtain
\begin{align*}
\Delta_1&=\left[(V'-rV)(m_1X_1-m_2X_2)
+rV^2\left(m_1X_1+m_2\frac{X_1X_2-X_3}{X_1}\right)\right]^2+4m_1m_2(V'-rV)^2(X_1X_2-X_3)\\
 &+4(V'-rV)rV^2m_2\left[ m_1(X_1X_2-X_3)+m_2X_2\frac{X_1X_2-X_3}{X_1} \right]>0.
\end{align*}
Hence $(P_2^t)$ has two eigenvalues.

Therefore, for all the above cases, the first eigenvalue $\lambda(t)$ of ($P^t_2$) exists.
Moreover, from the definition \eqref{two-weig-0-x1-defi-V} of $V$ and the equation \eqref{two-weig-0-x1-equi-cont}, it follows that $\lambda''(t)$ exists a.e. on $(0,x_1].$

For $t\in(x_1,x_2]$, equation \eqref{two-weig-equiv} takes the form
\begin{equation}\label{two-weig-x1-x2-w}
r[\alpha_2L(t)+\beta_2 R(t)][\mu_2 L(t)+\gamma_2 R(t)]=\alpha_3.
\end{equation}
Again applying  \eqref{LR} leads to
\begin{equation}\label{two-weig-x1-x2-LR-w}
rR^2(t)\left[\alpha_2\int^t_0\frac{{\rm d}s}{R^2(s)}+\beta_2\right]
\left[\mu_2\int^t_0\frac{{\rm d}s}{R^2(s)}+\gamma_2\right]=\alpha_3.
\end{equation}
By the definition \eqref{two-weig-0-x1-defi-V} of $V$, the equation \eqref{two-weig-x1-x2-LR-w} becomes
\begin{equation}\label{two-weig-x1-x2-equi-cont}
A_2 \lambda^2+B_2\lambda+C_2=0,
\end{equation}
where $$A_2=m_1m_2\left\{ -r(1-V)^2X_1X_2-X_3[rV(1-V)+V'] \right\},$$
$$ B_2=r\left[ m_1X_1+V^2m_2\frac{X_1X_2-X_3}{X_1}  \right]-(rV-V')(m_1X_1+m_2X_2),\ C_2=rV-V'.$$

We now prove the existence of the first eigenvalue of $(P_2^t)$. We divide three cases.

\noindent Case $1$. $A_2=0$. By \eqref{equ-condi-exis}, we have $A_2-B_2\neq0$, then $(P_2^t)$ has an eigenvalue.

\noindent Case $2$.  $A_2>0$.

\noindent If $C_2<0$, then $(P_2^t)$ has two eigenvalues from the property of quadratic function.

\noindent If $C_2>0$, from the fact that $X_1X_2>X_3>0$, it follows that the discriminant
\begin{align*}
\Delta_2=B_2^2-4A_2C_2
&>\left[(rV-V')(m_2X_2-m_1X_1)-r\left( m_1X_1+V^2m_2\frac{X_1X_2-X_3}{X_1}  \right)\right]^2\\
&+4m_1m_2(rV-V')^2(X_1X_2-X_3)>0.
\end{align*}
Thus $(P_2^t)$ has two eigenvalues.

\noindent Case $3$. $A_2<0$.

\noindent If $C_2>0$, then $(P_2^t)$ has two eigenvalues by the property of quadratic function.

\noindent If $C_2<0$, using $X_1X_2>X_3>0$, we obtain
\begin{align*}
\Delta_2=B_2^2-4A_2C_2
&>\left[(V'-rV)(m_1X_1-m_2X_2)+r\left(m_1X_1+V^2m_2\frac{X_1X_2-X_3}{X_1}\right)\right]^2\\
 &+4(V'-rV)rX_2m_2^2V^2\frac{X_1X_2-X_3}{X_1}>0.
\end{align*}
Hence $(P_2^t)$ has two eigenvalues.

Therefore the first eigenvalue $\lambda(t)$ of ($P^t_2$) exists.
Moreover, combining \eqref{two-weig-0-x1-defi-V} and \eqref{two-weig-x1-x2-equi-cont} implies that $\lambda''(t)$ exists a.e. on $(x_1,x_2].$

For $t\in(x_2,1)$, the equation \eqref{two-weig-equiv} becomes \begin{equation}\label{two-weig-x2-1-w}
r[\alpha_3 L(t)+\beta_3 R(t)]R(t)=\alpha_3.
\end{equation}
Now use the relation
\begin{equation}\label{RL}
R(t)=L(t)\int^1_t\frac{ds}{L^2(s)}
\end{equation}
to rewrite \eqref{two-weig-x2-1-w} as
\begin{equation}\label{two-weig-x2-1-LR-w}
r\left[\alpha_3+\beta_3\int^1_t\frac{{\rm d}s}{L^2(s)}\right]L^2(t)\int^1_t\frac{{\rm d}s}{L^2(s)}=\alpha_3.
\end{equation}
Define
\begin{equation}\label{two-weig-x2-1-defi-U}
U(t)=\frac{L_2}{R_2}\int^1_t\frac{{\rm d}s}{L^2(s)}, \ t\in[x_2,1],
\end{equation}
then the equation \eqref{two-weig-x2-1-LR-w} reduces to
\begin{equation}\label{two-weig-x2-1-equi-cont}
A_3 \lambda^2+B_3\lambda+C_3=0,
\end{equation}
where $$A_3=m_1m_2X_3(U'+rU-rU^2),$$
$$ B_3=-(U'+rU)(m_1X_1+m_2X_2)+rU^2\left(m_2X_2+m_1\frac{X_1X_2-X_3}{X_2}\right),\ C_3=U'+rU.$$

We now prove the existence of the first eigenvalue of $(P_2^t)$. We divide three cases.

\noindent Case $1$. If $A_3=0$, i.e., $U'+rU-rU^2=0$, then $B_3=-rU^2m_1\frac{X_3}{X_2}<0$, and hence $(P_2^t)$ has an eigenvalue.

\noindent Case $2$. If $A_3>0$, then $C_3>0$. From the fact that $X_1X_2>X_3>0$, it follows that the discriminant
\begin{align*}
\Delta_3=B_3^2-4A_3C_3&=\left[(U'+rU-rU^2)(m_1X_1-m_2X_2)
+rU^2m_1\frac{X_3}{X_2}\right]^2\\
 &+4m_1m_2(U'+rU-rU^2)^2(X_1X_2-X_3)>0,
\end{align*}
Thus $(P_2^t)$ has two eigenvalues.

\noindent Case $3$. $A_3<0$.

\noindent If $C_3>0$, then $(P_2^t)$ has two eigenvalues from the property of quadratic function.

\noindent If $C_3<0$, using $X_1X_2>X_3>0$, we obtain
\begin{align*}
\Delta_3&=\left[-(U'+rU)(m_2X_2-m_1X_1)
+rU^2\left(m_2X_2+m_1\frac{X_1X_2-X_3}{X_2}\right)\right]^2+4m_1m_2(U'+rU)^2(X_1X_2-X_3)\\
 &-4(U'+rU)rU^2m_1\left[ m_2(X_1X_2-X_3)+m_1X_1\frac{X_1X_2-X_3}{X_2} \right]>0,
\end{align*}
Hence $(P_2^t)$ has two eigenvalues.

Therefore, the first eigenvalue $\lambda(t)$ of ($P^t_2$) exists.
Moreover, by the definition \eqref{two-weig-x2-1-defi-U} of $U$ and equation \eqref{two-weig-x2-1-equi-cont}, the derivative $\lambda''(t)$ also exists a.e. on $(x_2,1).$
\qed

The coefficients in the piecewise representations of $\varphi$ and $\psi$ satisfy a system of algebraic equations that can be solved by the spectral data. The following lemma provides a rigorous statement of this fact.

\begin{lem}\label{lem-coeff-lamda-two}
Let $\lambda(t)$ be the first eigenvalue of the perturbed problem ($P^t_2$). Suppose \eqref{equ-w-1} is left-definite. Then
the quantities
\begin{equation}\label{5011}
\alpha_j,\ 1\le j\le 3;\ \beta_1=0,\ \beta_2\frac{R_j}{L_j},\ \beta_3\frac{R_j}{L_j},\ j=2,\ 3
\end{equation}
and
\begin{equation}\label{5012}
\gamma_j,\ 1\le j\le 3;\ \mu_3=0,\ \mu_1\frac{L_j}{R_j},\ \mu_2\frac{L_j}{R_j},\ j=2,\ 3
\end{equation}
are functions determined by $\lambda(t)$, where $L(x)$ and $R(x)$ be defined as in \eqref{ini-sols}, and $\alpha_j,\ \beta_j,\ \gamma_j,\ \mu_j,\ 1\le j\le 3$ are coefficients defined in \eqref{two-dirac-experss-varphi} and \eqref{two-dirac-experss-psi}.
\end{lem}
\prf
One the one hand, since $\lambda_1$ is the first eigenvalue of the problem ($P_2$), we have $\alpha_3(\lambda_1)=0$, where $\alpha_3$ is given in \eqref{equ-alph3}.
On the other hand, since $\lambda(x_1)$ is the first eigenvalue of the problem
$$
-y''+[q-r\delta(x-x_1)]y=\lambda [m_1\delta(x-x_1)+m_2\delta(x-x_2)]y,\ y(0)=0=y(1),
$$
or equivalently, of the problem
$$
-y''+qy=\lambda [\widetilde{m}_1\delta(x-x_1)+m_2\delta(x-x_2)]y,\ y(0)=0=y(1)
$$
with $\widetilde{m}_1=m_1+r/\lambda(x_1)$, we must have $\alpha_3(\lambda(x_1))=0$ with $m_1$ replaced by $\widetilde{m}_1$, that is,
\begin{equation}\label{504}
1-\lambda(x_1)[\widetilde{m}_1L_1R_1+m_2L_2R_2]
+\lambda^2(x_1)\widetilde{m}_1m_2L_1R_2[L_2R_1-L_1R_2]=0.
\end{equation}
Since $\lambda(x_2)$ is the first eigenvalue of the problem
$$-y''+qy=\lambda [\widetilde{m}_2\delta(x-x_2)+m_1\delta(x-x_1)]y,\ y(0)=0=y(1),$$
using $\widetilde{m}_2=m_2+r/\lambda(x_2)$, we obtain
\begin{equation}\label{505}
1-\lambda(x_2)[m_1L_1R_1+\widetilde{m}_2L_2R_2]
+\lambda^2(x_2)m_1\widetilde{m}_2L_1R_2[L_2R_1-L_1R_2]=0.
\end{equation}

Introducing the variables $X_1,\ X_2,\ X_3$ defined in \eqref{two-var-x1-x2-x3},
the condition $\alpha_3(\lambda_1)=0$ together with equations \eqref{504} as well as \eqref{505} forms the linear system
\begin{equation}\label{three-equ}
\left\{\aligned
&\lambda_1m_1X_1+\lambda_1m_2X_2-\lambda^2_1m_1m_2X_3=1,\\
&\lambda(x_1)\widetilde{m}_1X_1+\lambda(x_1)m_2X_2-\lambda^2(x_1)\widetilde{m}_1m_2X_3=1, \\
&\lambda(x_2)m_1X_1+\lambda(x_2)\widetilde{m}_2X_2-\lambda^2(x_2)m_1\widetilde{m}_2X_3=1.
\endaligned\right.
\end{equation}
After simplification, the system becomes
\begin{equation}
\left\{\aligned
&\lambda_1m_1X_1+\lambda_1m_2X_2-\lambda^2_1m_1m_2X_3=1,\\
&-\frac{m_2 r}{m_1}X_2+\lambda(x_1)\widetilde{m}_1m_2\left[\lambda_1-\lambda(x_1)\right]X_3
=1-\frac{\lambda(x_1)\widetilde{m}_1}{\lambda_1m_1}, \\
&(*)X_3=(**),
\endaligned\right.
\end{equation}
where
$$(*)=m_1\left\{ \lambda(x_2)m_2[\lambda_1-\lambda(x_2)]+\lambda(x_1)m_1[\lambda_1-\lambda(x_1)]
+r[\lambda_1-\lambda(x_1)-\lambda(x_2)] \right\},$$
$$(**)=\frac{m_1[\lambda_1-\lambda(x_1)]+m_2[\lambda_1-\lambda(x_2)]-r}{\lambda_1 m_2}.$$
Because the three equations in \eqref{three-equ} are simultaneously satisfied, spectral theory  guarantees the existence of a solution.
Since the condition ${\bf (D_2)}$ implies $X_3\neq 0$, it follows that $m_1[\lambda_1-\lambda(x_1)]+m_2[\lambda_1-\lambda(x_2)]\neq r$, i.e., $(**)\neq 0$. Thus the coefficient matrix and the augmented matrix of the system both have rank three, and so the solution is unique.
Hence $X_1,\ X_2$ and $X_3$ are uniquely determined by $\lambda_1$, $\lambda(x_1)$ and $\lambda(x_2)$.
Consequently,  all the quantities $L_1R_1$, $L_2R_2$ and $ L^2_1R^2_2$ can be viewed as the functions of $\lambda_1$, $\lambda(x_1)$ and $\lambda(x_2).$
So, we may set
\begin{equation}\label{509}
L_1R_1=h_1, \ L_2R_2=h_2,\ L^2_1R^2_2=h_3,
\end{equation}
where $h_j=h_j(\lambda_1,\lambda(x_1), \lambda(x_2)),\ 1\leq j\leq 3$.
Substituting \eqref{509} into \eqref{euq-two-calu-alph-beta} shows that all coefficients
\begin{equation}
\alpha_j,\ 1\le j\le 3;\ \beta_1=0,\ \beta_2\frac{R_j}{L_j},\  \beta_3\frac{R_j}{L_j},\ j=2,\ 3
\end{equation}
are determined by $\lambda(t)$.

Combining \eqref{euq-two-calu-mu-gamma} with \eqref{509} indicates that
\begin{equation}
\gamma_j,\ 1\le j\le 3;\ \mu_3=0,\ \mu_1\frac{L_j}{R_j},\  \mu_2\frac{L_j}{R_j},\ j=2, 3
\end{equation}
are also determined by $\lambda(t)$.
\qed

Now we are ready to state the reconstruction result for $q$ in $(P_2)$, which recovers the potential in the same form as presented in Remark \ref{rem-sing-dira}.

\begin{thm}\label{thm-recove-two}
Under the conditions in Lemma \ref{lem-exis-prop-lamd1-two}, let $\lambda(t)$ be the first eigenvalue function of ($P_2$). Then $q$ can be uniquely determined by $\lambda(t)$. Moreover, $q$ can be recovered explicitly by
$$
q(t)=\left\{\sqrt{1/Y'(t)}\right\}''\Big/\!\sqrt{1/Y'(t)}, \ t\in(0,1),
$$
where
 $Y(t)=\left\{\aligned
& V(t),\ t\in(0,x_2],\\
& U(t),\ t\in[x_2,1),
\endaligned\right.$
and  $Y$ satisfies the following Cauchy problems on the three subintervals:
\begin{equation}\label{3212}
\left\{\aligned
& Y'(t)=r\left[c_jY+b_j\right]\left[l_jY+g_j\right],\ t\in I_j=[x_{j-1},x_j],\\
& Y(t_j)=1,
\endaligned\right.
\end{equation}
for $1\leq j\leq 3$ and $t_1=t_2=x_1,\ t_3=x_2$. Here the coefficients $c_j,\ b_j,\ l_j,\ g_j,\ 1\le j\le 3$ are functions of $t$ completely determined by the given spectral data $\lambda(t)$, and their explicit forms are provided in \eqref{0x1}, \eqref{x1x2} and \eqref{514} below, that is
\begin{equation}
\left.\begin{aligned}
&c_1=1,\ b_1=0,\ l_1=\frac{\mu_1 L_1}{\gamma_1 R_1},\ g_1=1;\\
&c_2=\frac{\alpha_2}{\alpha_3},\ b_2=\frac{\beta_2R_1}{\alpha_3L_1},\ l_2=\frac{\mu_2L_1}{\alpha_3R_1},\ g_2=\frac{\gamma_2}{\alpha_3};\\
&c_3=-1,\ b_3=0,\ l_3=\frac{R_2\beta_3}{\alpha_3L_2},\ g_3=1.
\end{aligned}\right.
\end{equation}
\end{thm}
\prf
For $t\in(0,x_1]$, substituting \eqref{two-weig-0-x1-defi-V} into \eqref{two-weig-0-x1-LR-w} shows that $V$ satisfies the Cauchy problem
\begin{equation}\label{0x1}
\left\{\aligned
& V'(t)=rV\left[1+\frac{\mu_1 L_1}{\gamma_1 R_1}V\right],\ t\in(0,x_1],\\
& V(x_1)=1.
\endaligned\right.
\end{equation}
By Lemma \ref{lem-coeff-lamda-two}, all coefficients in the differential equation \eqref{0x1} are uniquely determined by $\lambda(t)$. And hence, Lemma \ref{lem-exis-prop-lamd1-two}, together with the standard theory of ODEs, guarantees the existence and uniqueness of the solution $V$ for $t\in(0,x_1]$.

For $t\in(x_1,x_2]$, substituting \eqref{two-weig-0-x1-defi-V} into \eqref{two-weig-x1-x2-LR-w} yields that $V$ satisfies the Cauchy problem
\begin{equation}\label{x1x2}
\left\{\aligned
& \alpha_3V'(t)=r\left[\alpha_2V+\beta_2\frac{R_1}{L_1}\right]
\left[\gamma_2+\mu_2\frac{L_1}{R_1}V\right],\ t\in[x_1,x_2],\\
& V(x_1)=1.
\endaligned\right.
\end{equation}
Again by Lemma \ref{lem-coeff-lamda-two} and \ref{lem-exis-prop-lamd1-two}, the coefficients are  determined solely by  $\lambda(t)$, guaranteeing a unique solution $V$ for $t\in(x_1,x_2]$.

For $t\in(x_2,1)$, substituting \eqref{two-weig-x2-1-defi-U} into \eqref{two-weig-x2-1-LR-w} gives that $U$ satisfies
\begin{equation}\label{514}
\left\{\aligned
& U'(t)=-rU\left[1+\frac{R_2\beta_3}{\alpha_3L_2}U\right],\ t\in[x_2,1),\\
& U(x_2)=1.
\endaligned\right.
\end{equation}
Once more, Lemma \ref{lem-coeff-lamda-two} and \ref{lem-exis-prop-lamd1-two} ensure that the coefficients depend only on $\lambda(t)$, guaranteeing the existence and uniqueness of the solution $U$.

On each subinterval, the potential can be expressed as $q=L''/L$ or $q=R''/R$. Using the relations between $L,\ R$ and the auxiliary variables $V,\ U$ derived from \eqref{two-weig-0-x1-defi-V} and \eqref{two-weig-x2-1-defi-U}, we obtain
$$
q(t)=\left\{\sqrt{1/V'(t)}\right\}''\Big/\!\sqrt{1/V'(t)},\ t\in(0,x_2),
$$
and
$$
q(t)=\left\{\sqrt{1/U'(t)}\right\}''\Big/\!\sqrt{1/U'(t)},\ t\in(x_2,1).
$$
Because the coefficients of all differential problems are uniquely fixed by $\lambda(t)$, the resulting $q$ is uniquely determined. This completes the proof.
\qed

\begin{rem}
As stated in Remark \ref{rem-comp}, we compare Theorem \ref{thm-recove-two} with \cite[Theorem~4.1]{ZHQC2025}.

For $t\in(0,x_1)$, using formula \eqref{two-weig-0-x1-equi-cont}, we have
$$\frac{\partial \lambda(t,0)}{\partial r}=- \frac{\sqrt{(m_1X_1+m_2X_2)^2-4m_1m_2X_3}- (m_1X_1-m_2X_2)}{2m_1X_1\sqrt{(m_1X_1+m_2X_2)^2-4m_1m_2X_3}} \frac{R_1}{L_1}L^2.$$
Substituting the above expresstion into \eqref{eqn:reconstruct-formula} leads to the relation  $q=L''/L$.

For $t\in(x_1,x_2)$, by \eqref{two-weig-x1-x2-equi-cont}, we have
$$\frac{\partial \lambda(t,0)}{\partial r}=-\frac{ (xL+yR)^2   }{ 2m_1m_2X_3^2\frac{R_1}{L_1}\sqrt{(m_1X_1+m_2X_2)^2-4m_1m_2X_3}   },$$
where
\begin{align*}
&d=\sqrt{(m_1X_1+m_2X_2)^2-4m_1m_2X_3},\ e=(X_1X_2-X_3)(m_1X_1+m_2X_2-m_2\frac{X_3}{X_1}),\\ &f=(X_1X_2-X_3)[2m_1m_2X_3-(m_1X_1+m_2X_2)^2+m_2 \frac{X_3}{X_1}(m_1X_1+m_2X_2)],\\ &e'=(m_1X_1+m_2X_2)X_1X_2-m_1X_1X_3\\
&f'=(m_1X_1+m_2X_2)m_1X_1X_3-(m_1X_1+m_2X_2)^2X_1X_2 +2m_1m_2X_1X_2X_3,\\
&x=\frac{R_1}{L_1}\sqrt{ de +f  },\
y=\sqrt{ de'+ f'}.
\end{align*}
Again substituting the above result  into \eqref{eqn:reconstruct-formula}, we further derive that  $q=(xL+yR)''/(xL+yR)$.

For $t\in(x_2,1)$, based on \eqref{two-weig-x2-1-equi-cont}, we get
$$\frac{\partial \lambda(t,0)}{\partial r}=- \frac{\sqrt{(m_1X_1+m_2X_2)^2-4m_1m_2X_3}- (m_1X_1-m_2X_2)}{2m_1X_1\sqrt{(m_1X_1+m_2X_2)^2-4m_1m_2X_3}} \frac{L_2}{R_2}R^2,$$
and similarly we can deduce $q=R''/R$.
\end{rem}

Note that $w(x)=m_1\delta(x-x_1)+m_2\delta(x-x_2)$ here, so $w\equiv 0$ on all the subintervals $(0, x_1)$, $(x_1, x_2)$ and $(x_2, 1)$.
Then, substituting all the above three formulations into \eqref{eqn:reconstruct-formula} can recover $q$ on $(0,1)$ except both of the points $x_1$ and $x_2$.

Accordingly, Theorem \ref{thm-recove-two} can also be regarded as a generalization of \cite[Theorem~4.1]{ZHQC2025} to the framework involving the two-point Dirac weights.

\bigskip
\section{For the case of multi-point Dirac weight}
\medskip

In the final section, we will discuss the general case of multi-point Dirac weight, that is, the problem $(P_n)$ with the weight
$$
w(x)=\sum^n_{j=1}m_j\delta(x-x_j),\ m_j>0,\ 0=x_0<x_1<x_2\cdots<x_n<x_{n+1}=1,\ 3\leq n\in \mathbb{N}.
$$

The solution $\varphi(x,\lambda)$ and $\psi(x,\lambda)$ defined as in \eqref{ini-sols-lamb} can be chosen  as
$$
\left\{\aligned
&\varphi(x,\lambda)=\alpha_j L(x)+\beta_j R(x),\  &x\in I_j=[x_{j-1},x_j],\\
&\psi(x,\lambda)=\mu_j L(x)+\gamma_j R(x), &x\in I_j=[x_{j-1},x_j],
\endaligned\right.
$$
where $\alpha_j=\alpha_j(\lambda)$, $\beta_j=\beta_j(\lambda)$, $\mu_j=\mu_j(\lambda)$, $\gamma_j=\gamma_j(\lambda)$ for $j=1,2,\cdots n+1$ and
$$
\alpha_1=1, \ \beta_1=0;\ \mu_{n+1}=0,\ \gamma_{n+1}=1.
$$
We can know that
\begin{equation}\label{equ-wronsk-n}
W[\varphi,\psi](x)\equiv-(\alpha_j\gamma_j+\beta_j\mu_j)(\lambda),\ 1\le j\le n+1,\ x\in[0,1].
\end{equation}
Applying the continuity of $\varphi$ and $\psi$ and the jump conditions of their derivatives  at $x_j$
yields the recursive relations
\begin{equation}\label{equ-coe-recur-rela-n}
\left\{\aligned
& \alpha_{j+1}=\alpha_j-\lambda m_j(\alpha_j L_j+\beta_j R_j)R_j, \\
& \beta_{j+1}=\beta_j+\lambda m_jL_j(\alpha_j L_j+\beta_j R_j), \   1\le j\le n
\endaligned\right.
\end{equation}
where $L_j=L(x_j),\ R_j=R(x_j)$. Note that $\alpha_{n+1}(\lambda)$ is the characteristic function of $(P_n)$, that is $\lambda$ is an eigenvalue
of $(P_n)$ if and only if $\alpha_{n+1}(\lambda)=0$.

From Theorem \ref{thm-equi} and \eqref{equ-wronsk-n}, we know that $\widetilde{\lambda}(t)$ is an eigenvalue of $(P^t_n)$
for some $q$ if and only if
\begin{equation}\label{equ-n-char-funct-perted}
r\varphi(t,\widetilde{\lambda}(t))\psi(t,\widetilde{\lambda}(t))=\alpha_{n+1}(\widetilde{\lambda}(t)),\ t\in[0,1].
\end{equation}
\noindent
{\bf Question 1.}
What additional natural conditions guarantee the existence of the first eigenvalue of $(P^t_n)$? What properties does these eigenvalues have?

Assume that the first eigenvalue $\lambda(t)$ of $(P^t_n)$ exist for $t\in(0,1)$.
On the one hand, for $t\in (0,x_n]$, define
\begin{equation}\label{equ-n-weig-defi-V}
V_j(t)=\frac{R_{j}}{L_{j}}\int^t_0\frac{{\rm d}s}{R^2(s)},\ t\in(x_{j-1},x_j],\ j=1\cdots n,
\end{equation}
which, together with \eqref{LR}, tells us that, the equation \eqref{equ-n-char-funct-perted} means that $V_j$ satisfies the Cauchy problem
\begin{equation}\label{equ-cauchy-n-Vj}
\left\{\aligned
& \alpha_{n+1}V_j'(t)=r\left[\alpha_jV_j+\beta_j\frac{R_{j}}{L_{j}}\right]
\left[\gamma_j+\mu_j\frac{L_{j}}{R_{j}}V_j\right],\ t\in(x_{j-1},x_j],\ j=1\cdots n, \\
& V_j(x_j)=1.
\endaligned\right.
\end{equation}
where $\alpha_j=\alpha_j(\lambda(t)),\ \beta_j=\beta_j(\lambda(t))$,  $\mu_j=\mu_j(\lambda(t))$ and $\gamma_j=\gamma_j(\lambda(t))$.

For $t\in (x_n,1)$, define
\begin{equation}\label{equ-n-weig-defi-U}
U_n(t)=\frac{R_{n}}{L_{n}}\int^1_t\frac{{\rm d}s}{L^2(s)}.
\end{equation}
Then, by \eqref{equ-n-weig-defi-U} and \eqref{RL}, the equation \eqref{equ-n-char-funct-perted} implies that $U_n$ satisfies the Cauchy problem
\begin{equation}\label{equ-cauchy-n-Un}
\left\{\aligned
&U_n'(t)=-r\left[\frac{\beta_{n+1}L_n}{\alpha_{n+1}R_n}U_n+1\right],\ t\in(x_n,1),\\
& U_n(x_n)=1.
\endaligned\right.
\end{equation}

On the other hand, for any fixed $k\in\{1, 2, \cdots, n\}$, we introduce $\alpha^k_{j}$ and $\beta^k_{j}$ satisfying the new recurrence relation as follows:
\begin{equation}
\left\{\aligned
& \alpha^k_{j+1}=\alpha^k_j-\lambda \widehat{m}_j(\alpha^k_j L_j+\beta^k_j R_j)R_j, \\
& \beta^k_{j+1}=\beta^k_j+\lambda \widehat{m}_jL_j(\alpha^k_j L_j+\beta^k_j R_j),\ 1\leq j \leq n,
\endaligned\right.
\end{equation}
where
\begin{equation}
\widehat{m}_j=\left\{\aligned
& m_j,\ j\neq k, \\
& m_k+\frac{r}{\lambda(x_k)},\ j=k.
\endaligned\right.
\end{equation}
Then \(\lambda(x_k)\) is an eigenvalue of the perturbed problem with the perturbation located at the mass point \(x_k\), if and only if
\[
\alpha_{n+1}^k(\lambda(x_k))=0,\quad k=1,\dots,n.
\]
Setting \(\alpha_{n+1}^0(\lambda):=\alpha_{n+1}(\lambda)\), we obtain the algebraic system
\begin{equation}\label{equ-n-chara-functi}
\alpha_{n+1}^k(\lambda(x_k))=0,\quad k=0, 1, \dots,n.
\end{equation}

If all coefficients $\alpha_j$, $\beta_j$, $\mu_j$ and $\gamma_j$ appearing in \eqref{equ-cauchy-n-Vj} and \eqref{equ-cauchy-n-Un} (or equivalently, $L_j,\ R_j$ by recursive relations \eqref{equ-coe-recur-rela-n}), are uniquely determined by $\lambda(x_j)$ via \eqref{equ-n-chara-functi}, then the standard theory of ODEs guarantees the existence and uniqueness of the solution $V_j$ for $t\in(x_{j-1},x_j],\ j=1\cdots n$ and $U_n$ for $t\in(x_n,1)$.
On each subinterval, the potential can be expressed as $q=L''/L$ or $q=R''/R$. Using relations \eqref{LR}, \eqref{RL}, \eqref{equ-n-weig-defi-V} and \eqref{equ-n-weig-defi-U}, we obtain
$$
q(t)=\left\{\sqrt{1/V_j'(t)}\right\}''\Big/\!\sqrt{1/V_j'(t)},\ t\in(x_{j-1},x_j),\ j=1\cdots n,
$$
and
$$
q(t)=\left\{\sqrt{1/U_n'(t)}\right\}''\Big/\!\sqrt{1/U_n'(t)},\ t\in(x_n,1).
$$
Thus, reconstructing the potential only requires solving some algebraic equations such as
\eqref{equ-n-chara-functi}, which involve the $2n$ unknown $\{L_j,  R_j\}$ $(1\le j\le n)$, together with Cauchy problems such as \eqref{equ-cauchy-n-Vj} and \eqref{equ-cauchy-n-Un}.
For these unknowns, normalizing one of them (e.g., setting $L_1=1$) reduces the number of unknowns to $(2n-1)$, while the system provides only $n+1$ equations.
Thus, $(n-2)$ additional conditions are needed to determine the remaining variables.
This problem does not appear when $n\le 2$, but it is indeed the fundamental algebraic obstacle to unique reconstruction when $n\ge3$.

To end this paper, we introduce several questions, motivated by the above obstacle, which remain to be studied in the sequel.

\noindent
{\bf Question 2.}
Under what conditions can the above constants $\{L_j, R_j\}$ be uniquely recovered from $\lambda(t)$?

A potential remedy is to introduce two first eigenvalue functions, $\lambda(t, r_1)$ and $\lambda(t, r_2)$, corresponding to distinct coupling constants $r_1 \neq r_2$. This leads to the augmented system
\begin{equation}\label{605}
\alpha^k_{n+1}(r_1, \lambda(x_k)) = 0 \quad (0 \le k \le n), \qquad
\alpha^k_{n+1}(r_2, \lambda(x_k)) = 0 \quad (1 \le k \le n),
\end{equation}
which consists of $2(n+1)$ equations, exceeding the number of unknowns by two.
Accordingly, the following question arises naturally.

\noindent
{\bf Question 3.} What is the complete set of compatibility conditions that two functions $\lambda(t, r_1)$ and $\lambda(t, r_2)$ must satisfy in order to arise from the same potential $q$?

\bigskip
\bigskip

\section*{Acknowledgments}
\medskip

This research was partially supported by the National Natural Science Foundation of China [Grant numbers 12271299 and 12071254]  and Shandong Provincial Fund [ZR2024MA005]. The authors are grateful to Professor Xiaoping Yuan for his helpful discussions and guidance.

\end{document}